\documentclass[11pt]{article}
\usepackage{}
\usepackage[total={6in, 9in}]{geometry}
 \usepackage{amsfonts}
\usepackage{amsthm}
\usepackage{amssymb}
\usepackage{amsmath,amsfonts}
\usepackage{mathrsfs}
\usepackage{cases}
\usepackage{latexsym,bm}
\usepackage{indentfirst}
\usepackage{color}
\usepackage{ifpdf}
\usepackage{graphicx}
\usepackage{psfrag}
\usepackage{enumerate}
\usepackage{dsfont}
\usepackage[
pdfauthor={shan},
pdftitle={An Ore-type condition  for hamiltonicity in tough graphs and the extremal examples},
pdfstartview=XYZ,
bookmarks=true,
colorlinks=true,
linkcolor=blue,
urlcolor=blue,
citecolor=blue,
bookmarks=true,
linktocpage=true,
hyperindex=true
]{hyperref}

\newtheorem{THM}{\textbf{Theorem}}
\newtheorem{THMs}{\textbf{Theorem}}[section]

\newtheorem{LEM}{\textbf{Lemma}}
\newtheorem{CLA}{\textbf{Claim}}

\newcommand{\pf}{\noindent\textbf{Proof}.\quad}

\newcommand{\iC}{\overset{\leftharpoonup }{C}}

\newcommand{\oC}{\overset{\rightharpoonup }{C}}

\newcommand{\CC}{\mathcal{C}}

\def\dist{{\fam0 dist}}
\newcommand{\arxiv}[1]{\href{http://arxiv.org/abs/#1}{\texttt{arXiv:#1}}}

\begin{document}
\title{An Ore-type condition  for hamiltonicity in tough graphs and the extremal examples}
\author{Masahiro Sanka\thanks{Partially supported by JST Doctoral Program Student Support Project Grant Number JPMJST2123 and JST ERATO Grant Number JPMJER2301, Japan.}\\
	\medskip Keio University, 3-14-1 Hiyoshi, Kohoku-ku, Yokohama 223-8522, Japan \\
	\medskip 
	{\tt sankamasa@keio.jp}
	\\ Songling Shan\thanks{Partially supported by NSF grant DMS-2345869, USA.}\\ 
	\medskip Aunurm  University, Auburn, AL 36849, USA\\
	\medskip 
	{\tt szs0398@auburn.edu}
}
\date{\today}
\maketitle

\emph{\textbf{Abstract}.}
Let $G$ be a $t$-tough graph on $n\ge 3$ vertices for some $t>0$.  It was  shown by Bauer et al. in 1995 that if the minimum degree of $G$ is greater than $\frac{n}{t+1}-1$, then $G$ is hamiltonian. 
In terms of Ore-type hamiltonicity conditions, the problem was only studied when $t$ is between 1 and 2, and recently the author 
proved a general result. The result states that if the degree sum of any two nonadjacent vertices of $G$ is greater than $\frac{2n}{t+1}+t-2$, then $G$ is hamiltonian. 
It was conjectured in the same paper that the ``$+t$" in the bound $\frac{2n}{t+1}+t-2$ can be removed.  Here we confirm 
the conjecture.  The result 
generalizes the result by Bauer, Broersma, van den Heuvel, and Veldman. Furthermore, we characterize  all $t$-tough graphs  $G$ on $n\ge 3$ vertices for which $\sigma_2(G) = \frac{2n}{t+1}-2$ but
$G$ is non-hamiltonian.


\vspace{2mm}

\section{Introduction}

We consider only simple  graphs. 
Let $G$ be a graph.
Denote by $V(G)$ and  $E(G)$ the vertex set and edge set of $G$,
respectively. Let $v\in V(G)$, $S\subseteq V(G)$, and $H\subseteq G$. 
Then  $N_G(v)$   denotes the set of neighbors
of $v$ in $G$, $d_G(v):=|N_G(v)|$ is the degree of $v$ in $G$, 
and $\delta(G):=\min\{d_G(v): v\in V(G)\}$ is the minimum degree of  $G$. 
Define 
$\deg_G(v,H)=|N_G(v)\cap V(H)|$, $N_G(S)=(\bigcup_{x\in S}N_G(x))\setminus S$, 
and  we write $N_G(H)$ for $N_G(V(H))$. 
Let $N_H(v)=N_G(v)\cap V(H)$ 
and $N_H(S)=N_G(S)\cap V(H)$.  Again, we write  $N_H(R)$ for $N_H(V(R))$
for any subgraph $R$ of $G$. 
We use $G[S]$ and $G-S$ to denote the subgraphs of $G$ induced by  $S$ and $V(G)\setminus S$, respectively. 
For notational simplicity we write $G-x$ for $G-\{x\}$.
Let $V_1,
V_2\subseteq V(G)$ be two disjoint vertex sets. Then $E_G(V_1,V_2)$ is the set
of edges in $G$  with one endvertex  in $V_1$ and the other endvertex  in $V_2$. For two integers $a$ and $b$, let $[a,b]=\{i\in \mathbb{Z}\,:\,   a\le i \le b\}$.

Throughout this paper,  if not specified, 
we will assume $t$ to be a nonnegative real number. The number of components of  a graph $G$ is denoted by $c(G)$. 
The graph $G$ is said to be {\it $t$-tough\/} if $|S|\ge t\cdot
c(G-S)$ for each $S\subseteq V(G)$ with $c(G-S)\ge 2$. The {\it
	toughness $\tau(G)$\/} is the largest real number $t$ for which $G$ is
$t$-tough, or is  $\infty$ if $G$ is complete. This concept was introduced by Chv\'atal~\cite{chvatal-tough-c} in 1973.
It is  easy to see that  if $G$ has a hamiltonian cycle
then $G$ is 1-tough. Conversely,
Chv\'atal~\cite{chvatal-tough-c}
conjectured that
there exists a constant $t_0$ such that every
$t_0$-tough graph is hamiltonian.
Bauer, Broersma and Veldman~\cite{Tough-counterE} have constructed
$t$-tough graphs that are not hamiltonian for all $t < \frac{9}{4}$, so
$t_0$ must be at least $\frac{9}{4}$ if Chv\'atal's toughness conjecture is true.

Chv\'atal's toughness conjecture  has
been verified  for certain classes of graphs including 
planar graphs, claw-free graphs, co-comparability graphs, and
chordal graphs~\cite{Bauer2006}.   The classes 
also include $2K_2$-free graphs~\cite{2k2-tough,1706.09029,2103.06760}, 
 and $R$-free 
graphs for $R\in \{P_2\cup P_3, P_3\cup 2P_1, P_2\cup kP_1\}$~\cite{p2p3, 2107.08476, MR4455928,  2210.10408, 2303.09741}, where 
$k\ge 4$ is an integer.  
 In general, the conjecture is still wide open. 
In finding hamiltonian cycles in graphs, sufficient  conditions such as Dirac-type  and 
Ore-type conditions are the most classic ones.

\begin{THMs}[Dirac's Theorem~\cite{Dirac-theorem}]
	If $G$ is a graph on $n\ge 3$  vertices with $\delta(G) \ge \frac{n}{2}$, then $G$
	is hamiltonian. 
\end{THMs} 

Define $\sigma_2(G)=\min\{d_G(u)+d_G(v)\,:\,  \text{$u,v\in V(G)$ and they are nonadjacent}\}$ 
if $G$ is noncomplete, and define $\sigma_2(G)=\infty$ otherwise. Ore's Theorem, as a generalization of 
Dirac's Theorem, is stated below. 

\begin{THMs}[Ore's Theorem~\cite{Ore-Theorem}]\label{thm:ore}
	If $G$ is a graph on $n\ge 3$ vertices with $\sigma_2(G) \ge n$, then $G$
	is hamiltonian. 
\end{THMs}

Analogous to Dirac's Theorem, Bauer, Broersma, van den Heuvel, and Veldman~\cite{MR1336668} proved the following result 
by incorporating the toughness of the graph. 

\begin{THMs}[Bauer et al.~\cite{MR1336668}]\label{degree-tough}
	Let $G$ be a $t$-tough graph on $n\ge 3$ vertices.  If $\delta(G) > \frac{n}{t+1}-1$, then 
	$G$ is hamiltonian. 
\end{THMs}

A natural question here is whether we can  find an Ore-type condition  involving the toughness of $G$
that generalizes Theorem~\ref{degree-tough}.  Various theorems were proved prior to 
Theorem~\ref{degree-tough} by only taking $\tau(G)$ between 1 and 2~\cite{MR499116,BCL, MR1032634}. 
	Let $G$ be a $t$-tough graph on $n\ge 3$ vertices. 
The author showed in~\cite{2103.05146}  that if  $\sigma_2(G) > \frac{2n}{t+1}+t-2$, then $G$
is hamiltonian. It was also conjectured in~\cite{2103.05146} that $\sigma_2(G) > \frac{2n}{t+1}-2$  
is the right bound. In this paper, we confirm the conjecture. 
For any odd integer  $n\ge 3$,  the complete bipartite graph $G:=K_{\frac{n-1}{2}, \frac{n+1}{2}}$ is $\frac{n-1}{n+1}$-tough and 
satisfies  $\sigma_2(G)=n-1 =\frac{2n}{1+\frac{n-1}{n+1}}-2$.  However, $G$ is not hamiltonian. 
Thus,  
the degree sum condition  that $\sigma_2(G) > \frac{2n}{t+1}-2$    is  best possible for a  $t$-tough graph on at least three vertices
to be hamiltonian. 
In fact, for any odd integers $n\ge 3$, any graph from the family $\mathcal{H}=\{H_{\frac{n-1}{2}} + \overline{K}_{\frac{n+1}{2}}: \text{$H_{\frac{n-1}{2}}$  is any graph on $\frac{n-1}{2}$ vertices}\}$  
is an extremal graph, where ``$+$'' represents the join of two graphs.    
We also show that $\mathcal{H}$ is the only family of extremal graphs.

\begin{THM}\label{main}
	Let $G$ be a $t$-tough graph on $n\ge 3$ vertices.  Then the following statements hold. 
	\begin{enumerate}[(a)]
		\item  If $\sigma_2(G)> \frac{2n}{t+1}-2$, 
		then $G$ is hamiltonian.  
		\item If $\sigma_2(G)= \frac{2n}{t+1}-2$ and $G$ is not hamiltonian, 
		 then $G\in \mathcal{H}$. 
	\end{enumerate}

\end{THM}

%

The remainder of this paper is organized as follows: in Section 2, we introduce some notation and preliminary 
results,  and in Section 3, we prove Theorem~\ref{main}.

\section{Preliminary results}

Let $G$ be a graph and $\lambda$  be a positive integer. Following~\cite{MR696293}, a cycle $C$ of $G$ is 
a  \emph{$D_\lambda$-cycle} if every  component of $G-V(C)$ has order less than $\lambda$. 
Clearly, a $D_1$-cycle is just a hamiltonian cycle.
We denote by $c_\lambda(G)$  the number of components of $G$
with order at least $\lambda$, and write $c_1(G)$ just as $c(G)$. 
Two subgraphs $H_1$ and $H_2$ of $G$ are  \emph{remote} if they are disjoint  and there is no
edge of $G$ joining a vertex of $H_1$ with  a vertex of $H_2$. 
For a subgraph $H$ of $G$, let $d_G(H)=|N_G(H)|$ be the degree of  $H$ in $G$. 
We denote by $\delta_\lambda(G)$ the minimum
degree of a connected subgraph of order $\lambda$ in $G$. Again $\delta_1(G)$
is just $\delta(G)$.

\begin{LEM}[\cite{p2p3}]\label{lem:idependentset-size}
	Let $t>0$ and $G$ be a non-complete $n$-vertex $t$-tough  graph. Then $|W|\le \frac{n}{t+1}$ for every independent set $W$ in $G$. 
\end{LEM}

Denote by $\oC$ an orientation of $C$.
We assume that the orientation is clockwise throughout the rest of this paper. For $x\in V(C)$,
denote the immediate successor of $x$ on $\oC$ by $x^+$ and the immediate  predecessor of $x$ on $\oC$ by $x^-$.
We  use $N^+_{C}(x)$ to denote the set of immediate  predecessors for vertices from $N_C(x)$. 
For $u,v\in V(C)$, $u\oC v$  denotes the segment of $\oC$
starting at $u$, following $\oC$ in the orientation,  and ending at $v$.
Likewise, $u\iC v$ is the opposite segment of $\oC$ with endpoints as $u$
and $v$.  Let $\dist_{\oC}(u,v)$  denote the length of the path $u\oC v$.
For any vertex $u\in V(C)$ and any positive integer $k$, define 
$$
L_u^+(k)=\{v\in V(C)\,:\,   \dist_{\oC}(u,v) \in [1,k]\}
$$
to be the set of $k$ consecutive successors of $u$. Hereafter, all cycles under consideration are oriented, and we will  not 
distinguish between the notation $C$ and $\oC$.

The following lemma provides a way of extending a cycle $C$  provided that the vertices 
outside $C$ have many neighbors on $C$.  The proof follows from Lemma~\ref{lem:idependentset-size}
and is very similar to the proof of Lemma 10 in~\cite{p2p3}: if we assume instead that $C$ cannot be extended by including $x$,
then $N^+_C(x)\cup \{x\}$ is an independent set in $G$. 


\begin{LEM}\label{cycle-extendabilit y}
	Let $t>0$  and $G$ be an $n$-vertex  $t$-tough graph, 
	and let $C$ be a non-hamiltonian cycle of $G$. 
	If $x\in V(G)\setminus V(C)$ satisfies  $\deg_G(x, C)>\frac{n}{t+1}-1$,  then $G$ has a cycle $C'$ such that $V(C')=V(C)\cup \{x\}$. 	
\end{LEM}


A path $P$ connecting two vertices $u$ and $v$ is called 
a {\it $(u,v)$-path}, and we write $uPv$ or $vPu$ in order to specify the two endvertices of 
$P$. Let $uPv$ and $xQy$ be two paths. If $vx$ is an edge, 
we write $uPvxQy$ as
the concatenation of $P$ and $Q$ through the edge $vx$.

For an integer $\lambda \ge 1$, if a graph $G$ contains a $D_{\lambda+1}$-cycle $C$  but no $D_\lambda$-cycle,
then  $V(G) \setminus V(C) \ne \emptyset$. Furthermore,  $G-V(C)$ has a component  of order $\lambda$. 
The result below with  $d_G(H)$ replaced by $\delta_\lambda(G)$ and $H$ replaced by 
any component of $G-V(C)$ with order $\lambda$ was proved in~\cite[Corollary 7(a)]{MR1336668}. 

\begin{LEM}[\cite{2103.05146}]\label{lem:dominating-cycle}
	Let $G$ be a $t$-tough $2$-connected graph of order $n$.
	Suppose $G$ has a $D_{s+1}$-cycle   but no $D_s$-cycle for some integer $s\ge 1$. Let $C$ be a $D_{s+1}$-cycle of $G$ such that $C$ minimizes $c_p(G-V(C))$ prior to minimizing $c_q(G-V(C))$ for any $p,q\in [1,s]$ with $p>q$.  Then $n\ge (t+|V(H)|) (d_G(H)+1)$
	for any component $H$ of $G-V(C)$. 
\end{LEM}

The lemma below is the key to get rid of the  ``$+t$" in the  lower bound $\frac{2n}{t+1}+t-2$ on $\sigma_2(G)$
for guaranteeing the existence of a hamiltonian cycle~\cite{2103.05146}. 

\begin{LEM}\label{lem:path-degree}
	Let $G$ be a $t$-tough $2$-connected graph of order $n$. 
	Suppose that $G$ has a $D_{\lambda+1}$-cycle   but no $D_\lambda$-cycle for some integer $\lambda \ge 1$. Let $C$
	be a cycle of $G$.  
Then $G-V(C)$ has a component $H$ with order at least $\lambda$ such that 
		 $\deg_G(x,C) \le  \frac{n}{t+1}-\lambda$
		for some $x\in V(H)$. 
	\end{LEM}

\pf  Since $G$ has no $D_\lambda$-cycle, it is clear that $G-V(C)$ has a  component of order at least $\lambda$. 
We suppose to the contrary that for each  component $H$  with order at least $\lambda$ of $G-V(C)$ and each $x\in V(H)$, we have $\deg_G(x,C)  >  \frac{n}{t+1}-\lambda$.    Among all cycles $C'$
of $G$ that satisfy the two conditions below, we  may assume that $C$
is one that minimizes $c_p(G-V(C))$ prior to minimizing $c_q(G-V(C))$ for any $p \ge \lambda$ and any $q$ with $q<p$. 
\begin{enumerate}[(1)]
	\item each component of $G-V(C)$ either has order at most $\lambda-1$, or  \label{C1}
	\item  the component $H$ has  order at least $\lambda$  such that  for each $x\in V(H)$, we have $\deg_G(x,C)  >  \frac{n}{t+1}-\lambda$.  \label{C2}
\end{enumerate}

We take a component $H$ with order at least $\lambda$ and assume that $N_{C}(H)$ has size $k$
for some integer $k\ge 2$, and that the $k$ neighbors are $v_1, \ldots, v_k$ and appear in the same order along $\oC$. 
Note that $k> \frac{n}{t+1}-\lambda$ by our assumption. 
For each $i\in [1,k]$, 
and each $v\in V(v_i^+\oC v_{i+1}^-)$,  where $v_{k+1}:=v_1$, we let $\CC(v)$ be the set of components of $G-V(C)$
that have a vertex joining to $v$ by an edge in $G$.  As $N_C(H)\cap V(v_i^+\oC v_{i+1}^-)=\emptyset$, we have 
$ H \notin \CC(v)$. 
Let $w_i^*\in V(v_i^+\oC v_{i+1}^-)$ be the vertex with $\dist_{\oC}(v_i,w_i^*)$ minimum  
such that 
$$
\sum\limits_{D\in \bigcup\limits_{v\in V(v_i^+\oC w_i^*)} \CC(v)}|V(D)|+ |V(v_i^+\oC w_i^*)| \ge \lambda. 
$$
If such a vertex $w_i^*$ exists, 
let $L_{v_i}^*(\lambda)$  be the union of the vertex set $V(v_i^+\oC w_i^*)$
and all those vertex sets of graphs in $\bigcup\limits_{v\in V(v_i^+\oC w_i^*)} \CC(v)$;
if such a vertex $w_i^*$ does not exist, let $L_{v_i}^*(\lambda)=L_{v_i}^+(\lambda)$. 
Note that when $w_i^*$ exists, by its definition, $w_i^*\in V(v_i^+\oC v_{i+1}^-)$. Thus  $V(v_i^+\oC w_i^*)\cap V(v_j^+\oC w_j^*) =\emptyset$ if both $w_i^*$ and $w_j^*$ exist for distinct $i,j\in [1,k]$.

We will show that we can make the following assumptions:
\begin{enumerate}[(a)]
	\item  If for some $i \in [1,k]$,  it holds that $L_{v_i}^*(\lambda)=L_{v_i}^+(\lambda)$, then $\dist_{\oC}(v_i, v_j) \ge \lambda+1$ for any  $j\in[1,k]$ with $j\ne i$.
	Thus the vertex  $w_i^*$ exists for each $i\in [1,k]$. 
	
	\item $G[L_{v_i}^*(\lambda)]$ and $G[L_{v_j}^*(\lambda)]$ are pairwise remote for any distinct $i,j\in [1,k]$. 
\end{enumerate}

With Assumptions (a) and (b),  we can reach a contradiction as follows:  note that $G[L_{v_i}^*(\lambda)]$ and $G[L_{v_j}^*(\lambda)]$ are 
remote for any distinct $i,j\in [1,k]$
and $H$  and $G[L_{v_i}^*(\lambda)]$ are  remote for any $i\in [1,k]$.  Let $S=V(G)\setminus \left( (\bigcup_{i=1}^k L_{v_i}^*(\lambda)) \cup V(H)\right)$.  
Then $|S| \le n-(k+1) \lambda$ and $c(G-S) = k+1$. As $G$ is $t$-tough, we get 
$$
n-(k+1) \lambda \ge |S| \ge t\cdot c(G-S) = t(k+1), 
$$
giving $ k \le  \frac{n}{t+\lambda}-1$.  Since $n\ge (\lambda+t)(2t+1)$ by Lemma~\ref{lem:dominating-cycle} ($G$ has a $D_{\lambda+1}$-cycle $C'$ such that $G-V(C')$ has a component  $H'$ of order $\lambda$, and $d_G(H') \ge 2t$ by $G$ being $t$-tough), we get    
\begin{eqnarray*}
	\frac{n}{t+1}-\lambda -\left(\frac{n}{t+\lambda}-1 \right) &=&\frac{(\lambda-1)(n-(t+1)(t+\lambda))}{(t+1)(t+\lambda)} \ge 0, 
\end{eqnarray*}
and so 
$k \le  \frac{n}{t+\lambda}-1 \le \frac{n}{t+1}-\lambda$. This gives a contradiction to $k> \frac{n}{t+1}-\lambda$. 
Thus  we are only left to show Assumptions (a) and  (b). 
We show that if any one of the assumptions is violated, then we can 
 decrease  $c_p(G-V(C))$ for some  $p\ge \lambda$.  

For Assumption (a), 
if $L_{v_i}^*(\lambda)=L_{v_i}^+(\lambda)$ for some $i\in [1,k]$
but $\dist_{\oC}(v_i,v_{j})  \le \lambda$ for some $v_j\in N_{C}(H)$ with $j\ne i$,  
then there must exist two consecutive indices $i,j\in [1,k]$ such that $\dist_{\oC}(v_i,v_{j})  \le \lambda$. 
Thus we may just assume $j=i+1$, where the index is taken modulo $k$. 
Let $v_i^*, v_{i+1}^*\in V(H)$ such that $v_iv_i^*, v_{i+1}v_{i+1}^*\in E(G)$, and let $P$ be a $(v_i^*, v_{i+1}^*)$-path in $H$. 
Let 
$C_1=v_i\iC v_{i+1}v^*_{i+1}Pv_i^*v_i$.

Note that every component of $G-V(C)$ not having any vertex joining to a vertex from $v_i^+\oC v_{i+1}^-$ in $G$ is still a component 
of $G-V(C_1)$. Those components automatically satisfy Conditions~\eqref{C1} and~\eqref{C2} as listed in the beginning of this proof.  Vertices  in $v_i^+\oC v_{i+1}^-$  are contained in a distinct component of $G-V(C_1)$, and the component has order at most $\lambda-1$ 
by the assumption that $L_{v_i}^*(\lambda)=L_{v_i}^+(\lambda)$ and $\dist_{\oC}(v_i,v_{i+1})  \le \lambda$. 
Finally, as any vertex from each component of $H-V(v^*_{i+1}Pv_i^*)$  is not adjacent in $G$  to any vertex from $v_i^+\oC v_{i+1}^-$, we know that components of $H-V(v^*_{i+1}Pv_i^*)$ are components of $G-V(C_1)$, and that 
 $\deg_G(w, C_1) > \frac{n}{t+1}-\lambda$  for any $w\in V(H-V(v^*_{i+1}Pv_i^*))$.  
  Hence each 
component of $G-V(C_1)$ 
either has order at most $\lambda-1$  or 
is a  component of order at least $\lambda$ such that each vertex from the component has in $G$
more than  $\frac{n}{t+1}-\lambda$ neighbors on $C_1$. 
However, $c_{|V(H)|}(G-V(C_1)) <c_{|V(H)|}(G-V(C))$ and $c_q(G-V(C_1))=c_q(G-V(C))$
for any $q >|V(H)|$, contradicting the choice of $C$. 
Therefore we have Assumption  (a),  which implies that the vertex  $w_i^*$ exists for each $i\in [1,k]$. 

For Assumption (b),  suppose it is false. Then there exist distinct  $i,j\in [1,k]$ such that $G[L_{v_i}^*(\lambda)]$ and $G[L_{v_j}^*(\lambda)]$ 
are not remote. By the definition of remote subgraphs, we have  either  $L_{v_i}^*(\lambda) \cap L_{v_j}^*(\lambda) \ne \emptyset$
or $L_{v_i}^*(\lambda) \cap L_{v_j}^*(\lambda) = \emptyset$ but $E_G(L_{v_i}^*(\lambda), L_{v_j}^*(\lambda)) \ne \emptyset$. 
In order to achieve a contradiction, 
we first show the following general claim, call it Claim ($*$). 

Claim ($*$):  For any $r \in [1,\dist_{\oC}(v_i, w_i^*)]$ and $s\in [1,\dist_{\oC}(v_j, w_j^*)]$, if $L_{v_i}^*(r) \cap L_{v_j}^*(s) = \emptyset$, then 
$E_G(L_{v_i}^*(r), L_{v_j}^*(s)) = \emptyset$.  
Suppose otherwise that $E_G(L_{v_i}^*(r), L_{v_j}^*(s))  \ne  \emptyset$.  
Since there is no edge of $G$ connecting  any two components of $G-V(C)$,  $E_G(L_{v_i}^*(r), L_{v_j}^*(s)) \ne \emptyset$ implies that there exist $y\in V(v_i^+ \oC w_i^*)\cap L_{v_i}^*(r) $ and $z\in V(v_j^+ \oC w_j^*)\cap L_{v_j}^*(s) $ 
such that $yz \in E(G)$. 
We choose  $y\in  V(v_i^+ \oC w_i^*) \cap L_{v_i}^*(r) $ with $\dist_{\oC}(v_i,y)$ minimum
and $z\in V(v_j^+ \oC w_j^*) \cap L_{v_j}^*(s) $  with $\dist_{\oC}(v_j,z)$ minimum such that 
$yz\in E(G)$. By this choice of $y$ and $z$,  it follows that 
$E_G(V(v_i^+ \oC y^-), V(v_j^+ \oC z^-))=\emptyset$. Let $v_i^*, v_{j}^*\in V(H)$ such that $v_iv_i^*, v_{j}v_{j}^*\in E(G)$, $P$ be a $(v_i^*, v_j^*)$-path in $H$, 
and 
let $C_1=v_i\iC zy\oC v_jv_j^*Pv_i^*v_i$.  Note that no vertex  of $H$ is adjacent in $G$ to any  vertex of $v_i^+ \oC y^-$ or $v_j^+ \oC z^-$ by the fact that  $V(v_i^+ \oC y^- )\subseteq  V(v_i^+ \oC w_i^*)$  and  $V(v_j^+ \oC z^-) \subseteq  V(v_j^+ \oC w_j^*)$ and Assumption (a). 
By the assumption that $L_{v_i}^*(r) \cap L_{v_j}^*(s) = \emptyset$  and the definitions of $L_{v_i}^*(\lambda) $ and $L_{v_j}^*(\lambda)$, 
we know that $v_i^+ \oC y^-$  and $v_j^+ \oC z^-$ are respectively contained in distinct  components of $G-V(C_1)$ 
that each of order at most $\lambda -1$.   By the same reasoning as in proving Assumption (a), 
we know that each 
component of $G-V(C_1)$ 
has order at most $\lambda-1$  or 
is a  component such that each vertex from the component has in $G$
more than $\frac{n}{t+1}-\lambda$ neighbors on $C_1$. 
However, $c_{|V(H)|}(G-V(C_1)) <c_{|V(H)|}(G-V(C))$ and $c_q(G-V(C_1))=c_q(G-V(C))$
for any $q >|V(H)|$, contradicting the choice of $C$.   Thus Claim ($*$) holds.

Now let us get back to prove Assumption (b) by contradiction.   Assume first that 
$L_{v_i}^*(\lambda) \cap L_{v_j}^*(\lambda) \ne \emptyset$. 
Then there exist $v\in V(v_i^+ \oC w_i^*)$ and $u\in V(v_j^+ \oC w_j^*)$ such that 
$\CC(v)\cap \CC(u) \ne \emptyset$,  
we  then further choose $v$ closest  to $v_i$ and $u$ closest to $v_j$  along $\oC$ with the property.  Thus  for any $w_i \in V(v_i^+\oC v^-)$ and any 
$w_j\in V(v_j^+\oC u^-)$,  it holds that $\CC(w_i)\cap \CC(w_j)= \emptyset$.  
Let $D\in \CC(v)\cap \CC(u)$ and  $v', u'\in V(D)$ such that $vv', uu'\in E(G)$, and $P'$
be a $(v',u')$-path of $D$. Let $v_i^*, v_{j}^*\in V(H)$ such that $v_iv_i^*, v_{j}v_{j}^*\in E(G)$, and let $P$ be a $(v_i^*, v_j^*)$-path in $H$. 
Then  $C_1=v_iv_i^*Pv_j^*v_j\iC vv'P'u'u\oC v_i$ is a cycle.  
Since each of $V(v_i^+\oC v^-)$ and $V(v_j^+ \oC u^-)$ contains at most $\lambda-1$
vertices and they are proper subsets of  $ V(v_i^+ \oC w_i^*)$ and $  V(v_j^+ \oC w_j^*)$ respectively, by Assumption (a) above, we have $N_{C}(H)\cap (V(v_i^+\oC v^-) \cup V(v_j^+ \oC u^-))=\emptyset$.  
By the choices of $v$ and $u$ that for any $w_i \in V(v_i^+\oC v^-)$ and any 
$w_j\in V(v_j^+\oC u^-)$,  it holds that $\CC(w_i)\cap \CC(w_j)= \emptyset$, 
Claim ($*$) implies that 
the components of $G-V(C_1)$
that respectively contain $v_i^+\oC v^-$ and $v_j^+ \oC u^-$ are disjoint. 
Since $V(v_i^+ \oC v^-)$  is a proper subset of $ V(v_i^+ \oC w_i^*)$ and  $V(v_j^+ \oC u^-)$  is a proper subset of $ V(v_j^+ \oC w_j^*)$, it follows by the definitions of $L_{v_i}^*(\lambda)$ and $L_{v_j}^*(\lambda)$
that the components of $G-V(C_1)$
that respectively contain $v_i^+\oC v^-$ and $v_j^+ \oC u^-$  have order at most $\lambda -1$.
By the same reasoning as in proving Assumption (a), 
we know that each 
component of $G-V(C_1)$ 
has order at most $\lambda-1$  or 
is a component such that each vertex from the component has in $G$
more than $\frac{n}{t+1}-\lambda$ neighbors on $C_1$. 
However, $c_{|V(H)|}(G-V(C_1)) <c_{|V(H)|}(G-V(C))$ and $c_q(G-V(C_1))=c_q(G-V(C))$
for any $q >|V(H)|$, contradicting the choice of $C$. Thus we  must have $L_{v_i}^*(\lambda) \cap L_{v_j}^*(\lambda) = \emptyset$. 
Applying Claim ($*$) again with $r=s=\lambda$, we have $E_G(L_{v_i}^*(\lambda), L_{v_j}^*(\lambda)) = \emptyset$. 
Therefore, $G[L_{v_i}^*(\lambda)]$ and $G[L_{v_j}^*(\lambda)]$  are remote, contradicting our assumption. Thus Assumption (b) holds. 
\qed

\section{Proof of Theorem~\ref{main}}

We may assume that $G$ is not a complete graph. Thus $G$ is $\lceil 2 t\rceil$-connected as it is $t$-tough. Suppose to the contrary that $G$ is not hamiltonian.

\begin{CLA}\label{claim:2-con}
We may assume that  $G$ is 2-connected. 
\end{CLA}

\pf  Since $t>0$, $G$ is connected. Assume to the contrary that $G$ has a cutvertex $x$.  
By considering the degree sum of two vertices respectively from two components of $G-x$, 
we know that $\sigma_2(G) \le n-1$.  On the other hand, $G$ has a cutvertex implies 
 $t\le \frac{1}{2}$ and so $\sigma_2(G) \ge  \frac{2n}{t+1}-2\ge \frac{4n}{3}-2$.  If $\sigma_2(G) > \frac{4n}{3}-2$, 
then we get a contradiction to  $\sigma_2(G) \le n-1$ as $n\ge 3$. Thus 
we assume  $\sigma_2(G) = \frac{4n}{3}-2$, which contradicts $\sigma_2(G) \le n-1$ if  $n\ge 4$. 
Thus $n=3$ and so $G=P_3$, but this implies  $G\in \mathcal{H}$. 
\qed


Since  $G$
is $2$-connected,  Lemma~\ref{lem:dominating-cycle} implies 
\begin{equation}\label{eqn:lower-bound-on-n}
	n\ge (t+1) (\lceil 2 t\rceil+1). \nonumber
\end{equation}
Also as  $G$
is $2$-connected, 
$G$ contains cycles. 
Let  $\lambda \ge 0$  be the integer 
such that $G$ admits no $D_\lambda$-cycle but a $D_{\lambda+1}$-cycle. 
Then we choose   $C$   
to be a longest $D_{\lambda+1}$-cycle 
that minimizes $c_p(G-V(C))$ prior to minimizing $c_q(G-V(C))$ for any $p,q\in [1,\lambda]$ with $p>q$.
 As $G$ is not hamiltonian, 
we have 
$\lambda\ge 1$. 
Thus $V(G)\setminus V(C) \ne \emptyset$. 
Since $C$ is not a $D_\lambda$-cycle but a $D_{\lambda+1}$-cycle, $G-V(C)$
has a component  $H$ of order $\lambda$.   
Let 
\begin{equation}\label{eqn:def-of-W}
	W=N_C(H) \quad \text{and} \quad  \omega=|W|.\nonumber     
\end{equation}
Since $G$ is a connected $t$-tough  graph,  it follows that 
$
\omega \ge  \lceil 2 t\rceil. 
$
On the other hand,  Lemma~\ref{lem:dominating-cycle} implies that 
$
\omega \le \frac{n}{t+ \lambda}-1.
$

\begin{CLA}\label{claim:H-and-W-size}
	\begin{numcases}{}
		\lambda+ \omega  < \frac{n}{t+1} & \text{if $\lambda \ge 2$},   \nonumber \\ 
		\lambda+ \omega   \le \frac{n}{t+1} & \text{if $\lambda=1$}.  \nonumber
	\end{numcases}
\end{CLA}

\pf  If $\lambda=1$, then  the assertion holds by $\omega \le \frac{n}{t+ \lambda}-1$.
Thus we assume $\lambda\ge 2$ and 
assume to the contrary that $\lambda+ \omega   \ge  \frac{n}{t+1}$. 
Then we have $ n\le (\lambda+\omega)(t+1)$.  By Lemma~\ref{lem:dominating-cycle}, we have 
$n \ge (\lambda+t)(\omega+1)$.  Thus we have 
$$
(\lambda+t)(\omega+1) \le (\lambda+\omega)(t+1), 
$$
which implies $\lambda \omega + \lambda +t\omega +t \le \lambda t + \lambda +t \omega +\omega$ 
and so $(\lambda-1) \omega  \le (\lambda-1) t$. 
Since $\lambda \ge 2$, we get $\omega \le t$, a contradiction to $\omega \ge 2t$. 
Note that the argument above for $\lambda \ge 2$
holds for all components of $G-V(C)$ as Lemma~\ref{lem:dominating-cycle} 
holds for all components of $G-V(C)$. 
%
\qed 

\begin{CLA}\label{claim:one-component}
If 	$\sigma_2(G) \ge \frac{2n}{t+1}-2$, then $H$ is the only component of $G-V(C)$.
\end{CLA}

\pf Suppose $H^* \ne H$ is another component of $G-V(C)$. 
Then we have $d_G(x)+d_G(y) \ge \sigma_2(G)$ for any $x\in V(H)$ and $y\in V(H^*)$. 
Since $d_G(x) \le \lambda+\omega-1$ and $d_G(y) \le |V(H^*)|+|N_C(H^*)|-1$, 
 Claim~\ref{claim:H-and-W-size}
implies that $|V(H^*)|+|N_C(H^*)|>\sigma_2(G) -(\frac{n}{t+1}-1)+1\ge \frac{n}{t+1}$ if $\lambda \ge 2$.
Repeating exactly the same argument for $|V(H^*)|+|N_C(H^*)|$ as in the proof of Claim~\ref{claim:H-and-W-size}
leads to a contradiction.  

Thus we assume $\lambda=1$. We get the same contradiction as above if $\sigma_2(G) >\frac{2n}{t+1}-2$ or $\lambda+\omega <\frac{n}{t+1}$. 
Thus we have $\sigma_2(G) =\frac{2n}{t+1}-2$ and $\omega =\frac{n}{t+1}-1$ by Claim~\ref{claim:H-and-W-size}. 
Then  $H$ and $H^*$ each contains only one vertex, say  $x$ and $y$, respectively.   
We first claim that 
the vertex $y$ is adjacent in $G$ to 
at most one vertex from $W^+$.  For otherwise, suppose  there are distinct $u, v\in W^+$
such that $yu,yv\in E(G)$. 
Then $C^*=u^-\iC vyu \oC v^-xu^-$ is a $D_{\lambda+1}$-cycle of $G$ with $c_\lambda(G-V(C^*))<c_\lambda(G-V(C))$. This contradicts the choice of $C$.  

We then claim that the set $W^+$ is an independent set in $G$. For otherwise, suppose  there are distinct $u, v\in W^+$
such that $uv\in E(G)$.  
Then $C^*=u^-\iC vu \oC v^-xu^-$ is a $D_{\lambda+1}$-cycle of $G$ with $c_\lambda(G-V(C^*))<c_\lambda(G-V(C))$. This contradicts the choice of $C$.  

Now let $S=V(G)\setminus (W^+\cup V(H)\cup V(H^*))$. 
Then $c(G-S) \ge \omega+1$. However
$$
\frac{|S|}{c(G-S)} \le \frac{n-\omega-2}{\omega+1}=\frac{\frac{tn}{t+1}-1}{\frac{n}{t+1}}<t, 
$$
 a contradiction. 

Therefore,  $H$ is the only component of $G-V(C)$. 
\qed 

Since $H$ is the only component of $G-V(C)$,  every vertex $v\in V(C)\setminus W$ is only adjacent in $G$ to vertices on $C$. 
As vertices  from $V(C)\setminus W$ are nonadjacent in $G$ with vertices from $H$, we have 
\begin{equation}\label{eqn:degree-on-C}
	\deg_G(v,C) \ge \sigma_2(G)-(\omega+\lambda-1)\quad \text{for any $v\in V(C)\setminus W$}. 
\end{equation}
We construct the vertex sets $L_u^+$ for each $u\in W$ as follows:
\begin{numcases}{L_u^+=}
\{ v\in V(C): \dist_{\oC}(u,v)<\frac{n}{t+1}-\omega+1\} & \text{if $\sigma_2(G) = \frac{2n}{t+1}-2$};  \nonumber  \\
\{ v\in V(C): \dist_{\oC}(u,v) \le \frac{n}{t+1}-\omega+1\} & \text{if $\sigma_2(G) > \frac{2n}{t+1}-2$}.  \nonumber 
\end{numcases}

\begin{CLA}\label{claim:segments}
	\begin{enumerate}[(a)]
		\item If $\sigma_2(G) =\frac{2n}{t+1}-2$, then for any two distinct vertices $u,v\in W$,  we have $\dist_{\oC}(u,v) \ge  \frac{n}{t+1}-\omega+1$ and 
		$E_G(L_u^+, L_v^+)=\emptyset$.  
		\item If $\sigma_2(G) > \frac{2n}{t+1}-2$, then for any two distinct vertices $u,v\in W$,  we have $\dist_{\oC}(u,v) > \frac{n}{t+1}-\omega+1$ and 
		$E_G(L_u^+, L_v^+)=\emptyset$. 
	\end{enumerate}
\end{CLA}

\pf  We only show Claim~\ref{claim:segments}(a), as the proof for Claim~\ref{claim:segments}(b) follows the same argument by 
just using the strict  inequality. 
Let $u^*\in N_H(u), v^*\in N_H(v)$ and $P$ be a $(u^*,v^*)$-path of $H$. 
For the first part of the statement, it suffices to show that when we arrange the vertices of $W$
along $\oC$, for any two consecutive vertices $u$ and $v$ from the arrangement, we have 
$\dist_{\oC}(u,v) \ge  \frac{n}{t+1}-\omega+1$. Note that $V(u^+\oC v^-)\cap W=\emptyset$ for such pairs of $u$ and $v$. 
Assume to the contrary that  there are distinct $u,v\in W$ with $V(u^+\oC v^-)\cap W=\emptyset$  and $\dist_{\oC}(u,v)  < \frac{n}{t+1}-\omega+1$. 
Let $C^*=u\iC vv^*Pu^*u$. Since $H$ has order  $\lambda$ and $V(u^+\oC v^-)\cap W=\emptyset$,   $H-V(P)$ is a union of components of $G-V(C^*)$  that each is of order at most $\lambda-1$
and $u^+\oC v^-$ is a  component of $G-V(C^*)$
of order less than $\frac{n}{t+1}-\omega$ but at least $\lambda$ ($G$ has no $D_\lambda$-cycle). 
By~\eqref{eqn:degree-on-C}, for each vertex $x\in V(u^+\oC v^-)$, $\deg_G(x, C^*)>\sigma_2(G)-(\omega+\lambda-1)-(\frac{n}{t+1}-\omega-1) = \frac{n}{t+1}- \lambda$. 
This shows a contradiction to Lemma~\ref{lem:path-degree}. 

For the second part of the statement, 
we assume to the contrary that $E_G(L^+_u, L^+_v) \ne \emptyset$. 
Applying the first part,   we know that $\dist_{\oC}(u,v) \ge  \frac{n}{t+1}-\omega+1$
and $\dist_{\oC}(v,u) \ge  \frac{n}{t+1}-\omega+1$ (exchanging the role of $u$ and $v$). 
Thus $L^+_u\cap L^+_v=\emptyset$. 
We choose  $x\in L^+_u$ with $\dist_{\oC}(u,x)$ minimum
and $y\in L^+_v$ with $\dist_{\oC}(v,y)$ minimum such that 
$xy\in E(G)$. By this choice of $x$ and $y$,  it follows that 
$E_G(V(u^+ \oC x^-), V(v^+ \oC y^-))=\emptyset$.
Let $C^*=u\iC yx\oC vv^*Pu^*u$. Since $H$ is  of order $\lambda$ and 
no vertex of $H$ is adjacent in $G$ to any  vertex of $u^+ \oC x^-$ or $v^+ \oC y^-$ by the first part of the statement,  $H-V(P)$ is a union of components of $G-V(C^*)$  that each is of order at most $\lambda-1$. Also
$u^+ \oC x^-$  and $v^+ \oC y^-$ are components of $G-V(C^*)$ 
that each is of order less than $ \frac{n}{t+1}-\omega$ but at least one of them has order at least $\lambda$.  

Since $E_G(V(u^+ \oC x^-), V(v^+ \oC y^-))=\emptyset$, by~\eqref{eqn:degree-on-C}, for each vertex $w\in V(u^+ \oC x^-)\cup V(v^+ \oC y^-)$, $\deg_G(w, C^*)>\frac{n}{t+1}- \lambda$.
This shows a contradiction to Lemma~\ref{lem:path-degree}. 
\qed 

By Claim~\ref{claim:segments}, $G[L_u^+]$ and $G[L_v^+]$ are  remote for any two distinct $u,v\in W$. 
Furthermore, $H$  is remote with $G[L_u^+]$  for any $u\in W$.  Furthermore, 
we have $|L_u^+| \ge \frac{n}{t+1} -\omega $ if $\sigma_2(G) =\frac{2n}{t+1}-2$, and $|L_u^+| > \frac{n}{t+1} -\omega $ if $\sigma_2(G) >\frac{2n}{t+1}-2$. 
 Let $S=V(G)\setminus \left( ( \bigcup_{ u\in W} L_{u}^+) \cup V(H)\right)$.  
Then  $c(G-S) = \omega+1$ and 
	\begin{numcases}{}
	|S| < n-\omega \left(\frac{n}{t+1}-\omega \right)-\lambda& \text{if  $\sigma_2(G) > \frac{2n}{t+1}-2$},   \nonumber \\ 
	|S| \le n-\omega \left (\frac{n}{t+1}-\omega \right)-\lambda & \text{if  $\sigma_2(G)  =  \frac{2n}{t+1}-2$}.  \nonumber
\end{numcases}
 As $G$ is $t$-tough and so $|S| \ge tc(G-S)=t(\omega+1)$, we get 
 	\begin{numcases}{}
 	n >\omega \left(\frac{n}{t+1}-\omega +t\right)+\lambda+t & \text{if  $\sigma_2(G) > \frac{2n}{t+1}-2$},   \label{eqn:n-lower-bound2a} \nonumber\\ 
 	n\ge \omega \left(\frac{n}{t+1}-\omega +t\right)+\lambda+t & \text{if  $\sigma_2(G)  = \frac{2n}{t+1}-2$}. \label{eqn:n-lower-bound2b}   \nonumber \end{numcases}

\begin{CLA}\label{claim:degree-sum=}
It holds that $\sigma_2(G) = \frac{2n}{t+1}-2$, $\lambda =1$, and $\omega=\frac{n}{t+1}-1$. 
\end{CLA}

\pf  Note that we have $\omega \le \frac{n}{t+1}-\lambda \le \frac{n}{t+1}-1$ by Claim~\ref{claim:H-and-W-size}. 
Suppose to the contrary that $\sigma_2(G) >\frac{2n}{t+1}-2$,  $\lambda\ge 2$, or $\omega<\frac{n}{t+1}-1$. 
Now we have 
$$
n\ge \omega \left(\frac{n}{t+1}-\omega +t\right)+\lambda+t, 
$$
implying 
\begin{equation}\label{eqn:in-claim5}
\left(\frac{\omega}{t+1}-1  \right) n \le \omega (\omega -t) -\lambda-t. 
\end{equation}
The inequality~\eqref{eqn:in-claim5} cannot achieve equality when $\sigma_2(G)>\frac{2n}{t+1}-2$, 
since we have $n >\omega \left(\frac{n}{t+1}-\omega +t\right)+\lambda+t$  in the case. 
If $\omega <t+1$, then we have $\omega <2$ because $2t\le \omega <t+1$ implies $t<1$, a contradiction to Claim~\ref{claim:2-con}. 
Thus we have $\omega \ge t+1$, implying $\frac{\omega}{t+1}-1\ge 0$. Then by Claim~\ref{claim:H-and-W-size}, we have 
\begin{equation}\label{eqn:in-claim5-2}
	\left(\frac{\omega}{t+1}-1  \right) n \ge \left(\frac{\omega}{t+1}-1  \right) (\omega+\lambda)(t+1). 
\end{equation}
Note that if $\lambda \ge 2$ or $\omega <\frac{n}{t+1}-1$, then the inequality~\eqref{eqn:in-claim5-2} cannot achieve the equality. 
By the assumption for the contrary, at least one of the inequalities~\eqref{eqn:in-claim5} or~\eqref{eqn:in-claim5-2} cannot achieve the equality. 
Therefore, combining ~\eqref{eqn:in-claim5} and~\eqref{eqn:in-claim5-2}, we get 
$$
\omega (\omega -t) -\lambda-t>\left(\frac{\omega}{t+1}-1  \right) (\omega+\lambda)(t+1), 
$$
which implies 
\begin{eqnarray*}
\omega^2-\omega t-\lambda-t &>& \omega (\omega +\lambda) -(\omega +\lambda) (t+1) \\ 
&=& \omega^2 + \omega \lambda -\omega t - \omega -\lambda t -\lambda. 
\end{eqnarray*}
This gives $(\lambda-1)t> (\lambda-1) \omega$, leading to $0<0$ or $\omega <t$, a contradiction. 
\qed 

By Claim~\ref{claim:degree-sum=}, Theorem~\ref{main}(a) holds. In the rest of the proof, we show Theorem~\ref{main}(b). 
Let  $$W^*=W^+\cup V(H).$$  
Since $u^+\in L_u^+$ for each $u\in W$, 
 Claim~\ref{claim:segments} implies that $W^*$ is an independent set in $G$. 

\begin{CLA}\label{claim:adjacency}
	Every vertex in $V(G) \setminus W^*$ is adjacent in $G$
	to at least two vertices from $W^*$. 
\end{CLA}

\pf  Suppose to the contrary that there exists $x\in V(G) \setminus W^*$
such that $x$ is adjacent in $G$ to at most one vertex from $W^*$. 
Let  $S=V(G)\setminus (W^*\cup \{x\})$. 
Then $c(G-S) \ge \omega+1$. However
$$
\frac{|S|}{c(G-S)} \le \frac{n-\omega-2}{\omega+1}=\frac{\frac{tn}{t+1}-1}{\frac{n}{t+1}}<t, 
$$
a contradiction.  
\qed

\begin{CLA}\label{claim:W^+-degree}
	For every $v\in W^+$, we have $\deg_G(v,C)=\frac{n}{t+1}-1$ and $v$ is not adjacent in $G$ to any two consecutive vertices on $C$. 
\end{CLA}

\pf Since $\sigma_2(G)=\frac{2n}{t+1}-2$, we have $\deg_G(v,C) \ge \frac{n}{t+1}-1$ 	for every $v\in W^+$. 
As $W^*$ is an independent set in $G$, $v^+\not\in W^*$. 
By Claim~\ref{claim:adjacency},  $v^+$ is adjacent in $G$ to another vertex $u$ from $W^*$. 
If $\{u\}=V(H)$, then $C^*=v^-\iC v^+ u v^-$ is a $D_{\lambda+1}$-cycle of $G$
with $v$ being the only component of $G-V(C^*)$.  
Assume then that $u\in W^+$. Let $V(H)=\{x\}$. Then 
$C^*=v^+ u\oC v^-xu^-\iC v^+$ is a $D_{\lambda+1}$-cycle of $G$
with $v$ being the only component of $G-V(C^*)$.  

Again, since $G$
has no $D_{\lambda}$-cycle, it follows that $\deg_G(v,C^*)=\frac{n}{t+1}-1$ and $v$ is not adjacent in $G$ to any two consecutive vertices on $C^*$.   The claim follows as $\deg_G(v,C)=\deg_G(v,C^*)$ and two neighbors  of $v$
that are consecutive on $C$ will also be consecutive on $C^*$. 
\qed 

Our goal is to show that $N_C(W^+)=N_C(H)$. To do so, we investigate 
how vertices in $N_C(W^+)$ are located along $\oC$. We start with some definitions. 
A \emph{chord} of $C$ is an edge $uv$ with $u,v\in V(C)$ and $uv\not\in E(C)$. 
Two chords $ux$ and $vy$ of $C$ that do not share any endvertices are \emph{crossing}  
if the four vertices $u, x, v, y$ appear  along  $\oC$  in the order $u,v, x, y$ or $u, y, x, v$. 
For two distinct vertices $x,y\in N_C(W^+)$, we say $x$
and $y$ form a \emph{crossing}  if there exist distinct vertices $u, v\in W^+$
such that $ux$ and $vy$ are crossing chords of $C$.   

\begin{CLA}\label{claim:no-crossing}
	For  any two distinct $x,y\in N_C(W^+)$ with $xy\in E(C)$, it follows that $x$ and $y$ do not form any crossing. 
\end{CLA}

\pf Suppose to the contrary that for some  distinct $x,y\in N_C(W^+)$ with $xy\in E(C)$, the two vertices $x$
and $y$ form a crossing. Let $u,v\in W^+$ such that $yu,yv\in E(G)$. Assume, without loss of generality, that 
the four vertices $u,v,x,y$ appear in the order $u, v, x, y$ along $\oC$. Let $V(H)=\{w\}$. 
Then $ux\iC vy \oC u^- w v^-\iC u$ is a hamiltonian cycle of $G$, a contradiction to 
our assumption that 
$G$ is not hamiltonian. 
\qed 

\begin{CLA}\label{claim:gap-contains-W+}
	For any vertex $v\in W^+$ and any two distinct $x,y\in N_C(v)$,  $x\oC y$ contains a vertex from $W^+$. 
\end{CLA}

\pf  By Claim~\ref{claim:W^+-degree}, $x\oC y$ has at least three vertices. Suppose to the contrary that 
$x\oC y$ contains no vertex from $W^+$.  We furthermore choose $x$ and $y$ so that  $x\oC y$ contains no 
other vertex from $N_C(v)\setminus\{x,y\}$. Assume that the three vertices  $v, x, y$ appear in the order  
$v, x, y$ along $\oC$.  By Claim~\ref{claim:adjacency}, each internal vertex of $x\oC y$
is adjacent in $G$ to a vertex from $W^+$. Then by our selection of $x$ and $y$, we know that 
each internal vertex of $x\oC y$
is adjacent in $G$ to a vertex from $W^+\setminus\{v\}$. 
Applying  Claim~\ref{claim:no-crossing},  $x^+$ does not form a crossing with $x$, and so $x^+$  forms
 a crossing with $y$. Similarly,  $x^{++}$ does not form a crossing with $x^+$, and so forms
 a crossing with $y$.  Continuing this argument for all the internal vertices of $x^{++}\oC y$,
 we know that 
$y^-$ 
forms a crossing with $y$, a contradiction to Claim~\ref{claim:no-crossing}. 
\qed 

We assume that the $\omega$ neighbors  of  the vertex from $V(H)$ on $C$ are $v_1, \ldots, v_\omega$ and  they appear in the same order along $\oC$. 
For each $i\in [1, \omega]$, let $I_i=V(v_i\oC v_{i+1}) \setminus \{v_{i}\}$, where $v_{\omega+1}:=v_1$. 

\begin{CLA}\label{claim:gap}
	For every $v\in W^+$,  it holds that $N_C(v)=W$. 
\end{CLA}

\pf  Since $x\oC y$ contains a vertex
 from $W^+$ for any two distinct $x,y\in N_C(v)$ by Claim~\ref{claim:gap-contains-W+}, it follows that no $I_i$ can contain more than one vertex 
 from $N_C(v)$. Since $\deg_G(v,C)=\omega=|W^+|$ by Claim~\ref{claim:W^+-degree}  and $\{I_1,\ldots, I_\omega\}$ 
 is a partition of $V(C)$, the Pigeon-hole Principle implies that each $I_i$
 contains exactly one vertex from  $N_C(v)$. 
 

Assume to the contrary that $N_C(v) \ne W$.  Let $i\in [1,\omega]$ be the  index such that $\dist_{\oC}(v,v_i)$ is largest and 
 $vv_{i}\not\in E(G)$.   Note that the   index $i$ exists since $v^-\in W$ and  $vv^-\in E(G)$. In particular, 
 every vertex $u\in W\cap V(v_i^+ \oC v)$ is adjacent to $v$ by the choice of $i$. Let $z$ be the vertex 
 in $N_C(v)\cap I_{i-1}$.  We prove the four subclaims below.  Let  $V(H)=\{x\}$ in the rest arguments. 
 
 {\bf \noindent Claim A: $z=v_i^-$. }
 
 \proof[Proof of Claim A]
Suppose otherwise that $z\ne v_i^-$.  Then by Claim~\ref{claim:adjacency}, $z^+$
 is adjacent in $G$ to at least two vertices from $W^+$. By Claim~\ref{claim:no-crossing}, 
 $N_C(z^+)\cap W^+ \subseteq V(v_i^+ \oC v)$. Thus $z^+$ 
 is adjacent in $G$ to a vertex from $W^+\cap V(v_i^+ \oC v^-)$ as $z$ is the only neighbor of $v$ from $I_{i-1}$ in $G$. By repeating this procedure
 for all the vertices from $V(z^{++} \oC v_i^-)$ iteratively, we conclude that $v_i^-$
  is adjacent in $G$ to a vertex $u \in W^+\cap V(v_i^+ \oC v^-)$. 
  As $v_i^+v_i \in E(G)$ and $v_iv_i^- \in E(C)$, Claim~\ref{claim:W^+-degree} implies that $v_i^+$
  is not adjacent in $G$ to $v_i^-$. Thus we have $u\not\in \{v_i^+,v\}$.  However, since $u^-v\in E(G)$ by our choice of  the index $i$, 
  the cycle $xv^-\iC uv_i^- \iC vu^- \iC v_i x$ is in $G$ longer than $C$, a contradiction. Thus $z$ must be $v_i^-$. 
 \qed 
 
 {\bf \noindent Claim B: $v_{i+1} =v^-$. }
 
 \proof[Proof of Claim B]
 Suppose that $v_{i+1}  \ne v^-$.  Considering $v_{i+1}^+$ in the place of $v$ and applying Claim A to it, $v_{i+1}^+$ 
 must be adjacent to $v_i$ or $v_i^-$ (if  $v_{i+1}^+v_{i}\not\in E(G)$, then $i$ is the  index such that $\dist_{\oC}(v_{i+1}^+,v_i)$ is largest and 
 $v_{i+1}^+v_{i}\not\in E(G)$).  If $v_{i+1}^+ v_i^- \in E(G)$, then the cycle $xv^-\iC v_{i+1}^+ v_i^- \iC vv_{i+1} \iC v_i x$ 
 is in $G$ longer than $C$, a contradiction. Thus we have $v_{i+1}^+ v_i\in E(G)$. We consider the vertex $v^+$. 
 Since $W^-$ is independent in $G$ and $v$ is adjacent to  $v_i^- \in W^-$, we have $v\not\in W^-$. Thus $v^+ \not\in W$. 
 Then by Claim~\ref{claim:adjacency}, $v^+$ is adjacent in $G$ to a vertex $u\in W^+\setminus \{v\}$. However, 
 the cycle 
 \begin{numcases}{}
 xv_i\oC vv_i^- \iC uv^+\oC u^- x  & \text{if $u\in V(v^+ \oC v_{i-1}^+)$},   \nonumber \\
 xv_iv_{i+1}^+ \oC vv_i^- \iC v^+ v_i^+ \oC v_{i+1} x & \text{if $u=v_i^+$},  \nonumber \\  
 xv^-\iC uv^+ \oC v_i^- vu^-\iC v_i x  & \text{if $u\in V(v_{i+1}^+ \oC v^-)$}, \nonumber 
 \end{numcases}
 is in $G$ longer than $C$, a contradiction. 
 \qed

 {\bf \noindent Claim C: $ \omega \ge 4$. }
 
 \proof[Proof of Claim C]
 
 Since $G$ is 2-connected by Claim~\ref{claim:2-con}, suppose instead that $\omega \in [2,3]$. First, suppose $\omega =2$. 
 Since $v\in W^+$ is adjacent to a vertex in $W^-$ and $W^+$ is independent in $G$, we have $W^-\setminus W^+ \ne \emptyset$. 
 Also a vertex $u \in W^-\setminus W^+$ is adjacent to all vertices in $W^+$ by Claim~\ref{claim:adjacency}. Then 
$u^+\in W$ and so  $u^{++} \in W^+$ is adjacent to $u^+$ and $u$, contrary to Claim~\ref{claim:W^+-degree}. Next, suppose $\omega =3$. 
We let, without loss of generality,  $v=v_1^+$. Then Claim B implies $v_1^+v_3^-\in E(G)$. Note that since $W^+$ is independent in $G$, $v_3^-$
must not be $v_2^+$.  We also have $v_2^+ v_1 \not\in E(G)$, as otherwise $v_1 v_2^+ \oC v_3^-v_1^+ \oC v_2 xv_3 \oC v_1$ is in $G$
a cycle longer than $C$. Applying Claim A to $v_2^+$, we  get $v_2^+ v_1^- \in E(G)$. 
Similarly, $v_3^+ v_2 \not\in E(G)$, as otherwise $v_3^+ v_2 \iC v_1x v_3\iC v_2^+ v_1^- \iC v_3^+ $ is in $G$
a cycle longer than $C$. Applying Claim A to $v_3^+$, we  get $v_3^+ v_2^- \in E(G)$.  
Then as the degrees of all vertices from $W^+$ are of degree 3 in $G$, Claims~\ref{claim:one-component}, \ref{claim:adjacency}, and~\ref{claim:W^+-degree} 
imply that the graph $G$ is isomorphic to the Petersen graph. However, $3=\omega =\frac{10}{t+1}-1$ implies that $G$ is $\frac{3}{2}$-tough, contradicting 
that the toughness of the Petersen graph is at most $\frac{4}{3}$ (in the Petersen graph, deleting two independent vertices from one 5-cycle and another two independent vertices 
that are non-neighbors of the first two deleted vertices from the second disjoint 5-cycle gives three components).  Thus we have $\omega \ge 4$. 
 \qed

  {\bf \noindent Claim D: For every $j\in [1,\omega]$, $|I_j| \ne 3$.}
 
 \proof[Proof of Claim D] Suppose that $|I_j| =3$
 for some $j\in [1,\omega]$. Then we have $v_j^+ v_{j+1}^- \in E(C)$, which implies $N_C(v_j^+) \ne W$. Applying Claim A to $v_j^+$, 
 we get $v_j^+ v_{j-1}^- \in E(G)$. By symmetry of the orientation of $C$, we have $v_{j+1}^-v_{j+2}^+\in E(G)$. Also we have $\omega \ge 4$
 by Claim C, which implies $v_{j-1} \in V(v_{j+2}^+ \oC v_j^-)$.  Then the cycle $xv_{j-1} \oC v_j^+ v_{j-1}^- \iC v_{j+2}^+ v_{j+1}^- \oC v_{j+2} x$ is in $G$
 longer than $C$, a contradiction. 
 \qed 
 
 We now show a contradiction. The vertex  $v^+$ must not be in $W$  since $vv_i^- \in E(G)$  by Claim A and $W^-$ is independent in $G$. 
 Thus $v^+$ is adjacent in $G$ to a vertex $u\in W^+\setminus \{v\}$ by Claim~\ref{claim:adjacency}. If $u \ne v_i^+$, 
 then the cycle $xv_i \oC v v_i^- \iC  uv^+ \oC u^- x$ is in $G$ longer than $C$, a contradiction. Thus we have $u=v_i^+$. 
 We consider the cycle $C^*=v^+ \oC v_i^- v_i x v^- \iC v_i^+ v^+$ in $G$. Note that we have $V(C^*)=V(G) \setminus \{v\}$. Then since the length of $C^*$
 is equal to the length of $C$, we can apply Claim D to $C^*$. However, $v^-, x, v_i, v_i^-$ are four consecutive vertices on $C^*$ appearing  in the order  $v^-, x, v_i, v_i^-$ and $v^-, v_i^- \in N_{C^*}(v)$, 
 showing that $C^*$
 does not satisfy Claim D, a contradiction. This completes the proof of Claim~\ref{claim:gap}. 
\qed 

Claim~\ref{claim:gap} implies that $N_C(W^*)=W$.  
Thus  every vertex from $W^*$ is adjacent in $G$ to every vertex from $W$. 
Therefore $t\le \tau(G) \le \frac{|W|}{|W^*|}$ as $W^*$ is an independent set in $G$.  Consequently, $|W|\ge t|W^*| =\frac{tn}{t+1}$ and so $W=V(G)\setminus W^*$ 
by noticing  $|W^*|=\frac{n}{t+1}$.
Thus $G$ contains a spanning complete bipartite graph between $W^*$ and $W$. 
On the other hand, since $|W^+|=|W|=\frac{n}{t+1}-1$ and $V(G)=W^*\cup W=(W^+\cup V(H)) \cup W$, we know that  $2(\frac{n}{t+1}-1)+1=n$
and so $t=\frac{n-1}{n+1}$. Thus  $|W|=\frac{n-1}{2}$ and $|W^*|=\frac{n-1}{2}+1=\frac{n+1}{2}$. 
Therefore, $G\in \mathcal{H}$. 
The proof of Theorem~\ref{main} is now complete.  
\qed

\section*{Acknowledgements}

The authors are very grateful to the  anonymous referees  for their careful reading and valuable
comments, which greatly  improved this paper.


\end{document}